\documentclass[letter, 11pt]{article}
\usepackage[T1]{fontenc}
\usepackage{tgpagella}
\usepackage{graphicx,amsmath,amsfonts,amscd,amssymb,bm,url,color,latexsym}
\usepackage{fullpage}
\usepackage[small,bf]{caption}
\usepackage{subcaption}
\usepackage{bbm}
\usepackage{microtype}
\usepackage{algorithm}
\usepackage{algorithmicx}
\usepackage{algpseudocode} 
\usepackage{hyperref}
\usepackage{verbatim}
\usepackage{framed}
\usepackage{comment}
\usepackage{wrapfig}
\usepackage[rightcaption]{sidecap}

\allowdisplaybreaks  

\hypersetup{
    colorlinks=true,%
    citecolor=blue,%
    filecolor=blue,%
    linkcolor=blue,%
    urlcolor=blue
}
\usepackage[toc, page]{appendix}

\newtheorem{theorem}{Theorem}[section]

\newtheorem{definition}[theorem]{Definition}

\usepackage{tikz}
\usepackage{fge}

\renewcommand{\mathbf}{\boldsymbol}

\newcommand{\mb}{\mathbf}
\newcommand{\mc}{\mathcal}

\newcommand{\bb}{\mathbb}
 
\newcommand{\set}[1]{\left\{ #1 \right\}}

\newcommand{\reals}{\bb R}

\newcommand{\R}{\reals}
\newcommand{\Cp}{\bb C}
\newcommand{\Z}{\bb Z}

\newcommand{ \paren }[1]{ \left( #1 \right) }

\DeclareMathOperator{\rank}{rank}
\DeclareMathOperator{\trace}{tr}

\DeclareMathOperator{\grad}{grad}
\DeclareMathOperator{\Hess}{Hess}
\DeclareMathOperator{\mini}{minimize}
\DeclareMathOperator{\maxi}{maximize}
\DeclareMathOperator{\st}{subject\; to}

\newcommand{\e}{\mathrm{e}}
\newcommand{\im}{\mathrm{i}}

\newcommand{\wh}{\widehat}

\newcommand{\ol}{\overline}

\newcommand{\norm}[2]{\left\| #1 \right\|_{#2}}
\newcommand{\abs}[1]{\left| #1 \right|}
\newcommand{\innerprod}[2]{\left\langle #1,  #2 \right\rangle}

\numberwithin{equation}{section}

\def \endprf{\hfill {\vrule height6pt width6pt depth0pt}\medskip}

\title{When Are Nonconvex Problems Not Scary?}
\author{Ju Sun, Qing Qu, and John Wright \\
\texttt{\{js4038, qq2105, jw2966\}@columbia.edu} \\
  \vspace{-.1cm}\\
Department of Electrical Engineering, Columbia University, New York, USA
}
\date{October 20, 2015 \; Revised: \today}

\begin{document}
\maketitle

\vspace{-0.3in}
\begin{abstract}
In this note, we focus on smooth nonconvex optimization problems that obey: (1) all local minimizers are also global; and (2) around any saddle point or local maximizer, the objective has a negative directional curvature. Concrete applications such as dictionary learning, generalized phase retrieval, and orthogonal tensor decomposition are known to induce such structures. We describe a second-order trust-region algorithm that provably converges to a global minimizer efficiently, without special initializations. Finally we highlight alternatives, and open problems in this direction. 
\end{abstract}

\section{Introduction}
General nonconvex optimization problems (henceforth ``NCVX problems'' for brevity) are NP-hard, even the goal is computing only a local minimizer~\cite{murty1987some, bertsekas1999nonlinear}. In applied disciplines, however, NCVX problems abound, and heuristic algorithms such as gradient descent and alternating directions are often surprisingly effective. The ability of natural heuristics to find high-quality solutions for practical NCVX problems remains largely mysterious. 

In this note, we study a family of NCVX problems that {\em can} be solved efficiently. This family cuts across central tasks in signal processing and machine learning, such as complete (sparse) dictionary learning~\cite{sun2015complete_tr}, generalized phase retrieval~\cite{sun2015geometric}, orthogonal tensor decomposition~\cite{ge2015escaping}, and noisy phase synchronization and community detection~\cite{boumal2016nonconvex,bandeira2016low}. 

Natural optimization formulations for these distinct tasks are nonconvex; surprisingly they exhibit a common characteristic structure. In each case, the goal is to estimate or recover an object from observed data. Under certain technical hypotheses, {\em every local minimizer of the objective function exactly recovers the object of interest}. 

With this structure, the central issue is how to escape the saddle points and local maximizers.
\begin{SCfigure}[][!htbp]
\centering
\includegraphics[width = 0.55\linewidth]{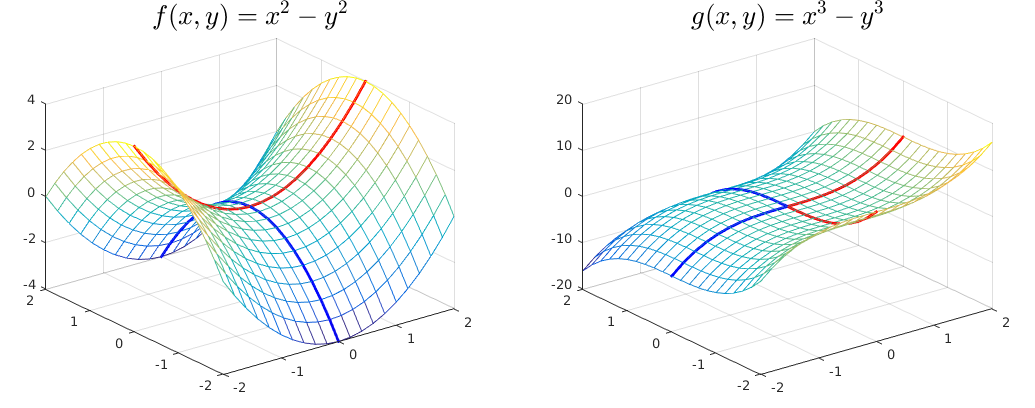}
\caption{Not all saddle points are ridable! Shown in the plot are functions $f(x, y) = x^2 - y^2$ (left) and $g(x, y) = x^3 - y^3$ (right). For $g$, both first- and second-order derivatives vanish at $(0, 0)$, producing a saddle that is induced by third-order derivatives. In both plots, red curves indicate local ascent directions and blue curves indicate local descent directions. }
\label{fig:illus_saddles}
\end{SCfigure}
Fortunately, for these problems, \emph{all saddle points and local maximizers are ``typical'' -- the associated Hessian matrix has at least one negative eigenvalue}. Geometrically, this means around any saddle point or local maximizer, the objective function has a negative curvature in a certain direction. Particularly, we call saddles of this type \emph{ridable saddles};\footnote{They are also called ``strict saddle" points in optimization literature, see, e.g., pp 38 of~\cite{reich2010nonlinear}; see also~\cite{ge2015escaping}.} the importance of this apparently extraneous restriction is illustrated in Figure~\ref{fig:illus_saddles}. Intuitively, at saddle points or local maximizers, in the direction of negative curvature the objective function is also locally descending. One can use this to design algorithms that escape from the saddle points and local maximizers concerned here. Indeed, consider a natural quadratic approximation to the objective $f$ around a saddle point $\mb x$: 
\begin{align*}
\wh{f}(\mb \delta; \mb x) =  f(\mb x) + \frac{1}{2} \mb \delta^* \nabla^2 f(\mb x) \mb \delta.
\end{align*}
When $\mb \delta$ is chosen to align with one eigenvector associated with a negative eigenvalue $\lambda_{\mathrm{neg}} [\nabla^2 f(\mb x)]$, it holds that 
\begin{align*}
\wh{f}(\mb \delta; \mb x) - f(\mb x) \le -\frac{1}{2} |\lambda_{\mathrm{neg}}| \norm{\mb \delta}{}^2.
\end{align*} 
Thus, minimizing $\wh{f}(\mb \delta; \mb x)$ locally provides a direction $\mb \delta_\star$ that tends to decrease the objective $f$, provided local approximation of $\wh{f}$ to $f$ is reasonably accurate.\footnote{For general saddles that seem to demand higher-order approximations, the computation may quickly become intractable. For example, third order saddle points seem to generally demand studying spectral property of three-way tensors, which entails NP-hard computational problems~\cite{hillar2013most}; see~\cite{anandkumar2016efficient} for a recent attempt in this line. } Based on this intuition, we derive an algorithmic framework that can exploit the second-order information to escape from saddle points and local maximizers and provably returns a global minimizer.

\section{Nonconvex Optimization with Ridable Saddles}
In this section, we present a more quantitative definition of the problem class we focus on and provide several concrete examples in this class.  

We are interested in optimization problem of the form: 
\begin{align} \label{eq:manifold_min}
\mini f(\mb x), \quad \st \quad \mb x \in \mc M. 
\end{align}
Here we assume $f$ is twice continuously differentiable, i.e., it has continuous first- and second-order derivatives, and $\mc M$ is a Riemannian manifold. Restricting $f$ to $\mc M$ and (with abuse of notation) writing the restricted function as $f$ also, one can effectively treat~\eqref{eq:manifold_min} as an unconstrained optimization on $\mc M$. We further use $\grad f(\mb x)$ and $\Hess f(\mb x)$ to denote the Riemannian gradient and Hessian of $f$ at point $\mb x$~\footnote{Detailed introduction to these quantities can be found in~\cite{absil2009}. We prefer to keep this at an intuitive level not to obscure the main ideas. }, which one can think of as Riemannian counterparts of Euclidean gradient and Hessian for functions, with the exception that $\grad f(\mb x)[\cdot]$ and $\Hess f(\mb x)[\cdot]$ only act on vectors in tangent space of $\mc M$ at $\mb x$, i.e., $T_{\mb x} \mc M$. 

\begin{definition}[($\alpha, \beta, \gamma, \delta$)-$\mc X$ functions]
A function $f: \mc M \mapsto \R$ is ($\alpha, \beta, \gamma, \delta$)-$\mc X$ ($\alpha$, $\beta$, $\gamma$, $\delta > 0$ ) if:\\
\indent 0) all local minimizers of $f$ are also global minimizers; \\
and $f$ is $(\alpha, \beta, \gamma, \delta)$-saddle,\footnote{See also strict-saddle function defined in~\cite{ge2015escaping}. When $\mc M$ is $\R^n$ or $\Cp^n$, the two definitions coincide. It is interesting to see if the two agree in general settings. Particularly, \cite{ge2015escaping} deals only with sets defined by equalities of the form $c_i(\mb x) = 0$ with differentiable function $c$, which excludes many manifolds of interest, such as symmetric positive definite matrices of a fixed dimension. See this page: \url{http://www.manopt.org/tutorial.html\#manifolds} for more examples. See also discussion in Introduction of this paper~\cite{absil2007trust} on relationship between manifold optimization and constrained optimization in the Euclidean space.} i.e., any point $\mb x \in \mc M$ obeys \textbf{at least one of the following}: ($T_{\mb x} \mc M$ is the tangent space of $\mc M$ at point $\mb x$)\\
\indent 1) [\textbf{Strong gradient}] $\norm{\grad f(\mb x)}{} \ge \beta$; \\
\indent 2) [\textbf{Negative curvature}] There exists $\mb v \in T_{\mb x} \mc M$ with $\norm{\mb v}{} = 1$ such that $\innerprod{\Hess f(\mb x) [\mb v]}{\mb v} \le -\alpha$; \\
\indent 3) [\textbf{Strong convexity around minimizers}] There exists a local minimizer $\mb x_\star$ such that $\norm{\mb x - \mb x_\star}{} \le \delta$, and for all $\mb y \in \mc M$ that is in $2\delta$ neighborhood of $\mb x_\star$, $\innerprod{\Hess f(\mb y)[\mb v]}{\mb v} \ge \gamma$ for any $\mb v \in T_{\mb y} \mc M$ with $\norm{\mb v}{} =1$, i.e., the function $f$ is $\gamma$-strongly convex in $2\delta$ neighborhood of $\mb x_\star$. \footnote{Strong convexity is required for the sake of deriving concrete convergence rate near the minimizers. Less stringent conditions might already be sufficient, depending on the specific problems and target computational guarantees. Also, it is possible to modify the trust-region methods (described later) to take advantage of fine problem structure; see, e.g.,\cite{sun2015geometric}.}
\end{definition}

In words, the function has no spurious local minimizers. Moreover, each point on the manifold $\mc M$ either has strong Riemannian gradient, or has Riemannian Hessian with at least one strictly negative eigenvalue, or lives in a small neighborhood of a local minimizer, such that the function is locally strongly convex. We remark in passing that requiring a function to be ridable may appear far too restrictive than it actually is. Indeed, one of the central results in Morse theory implies that a generic smooth function is ridable.  

In this note, we deal exclusively with minimizing $\mc X$ functions. \footnote{Though if the target is to compute any local minimizer of a ridable function, our method also applies. } These functions indeed appear in natural nonconvex formulations of important practical problems.  

\begin{SCfigure}[][!htbp]
\centering
\includegraphics[width = 0.35\linewidth]{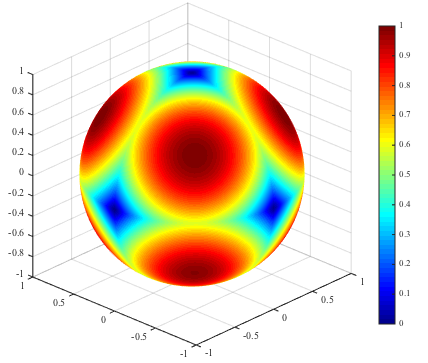}
\includegraphics[width = 0.35\linewidth]{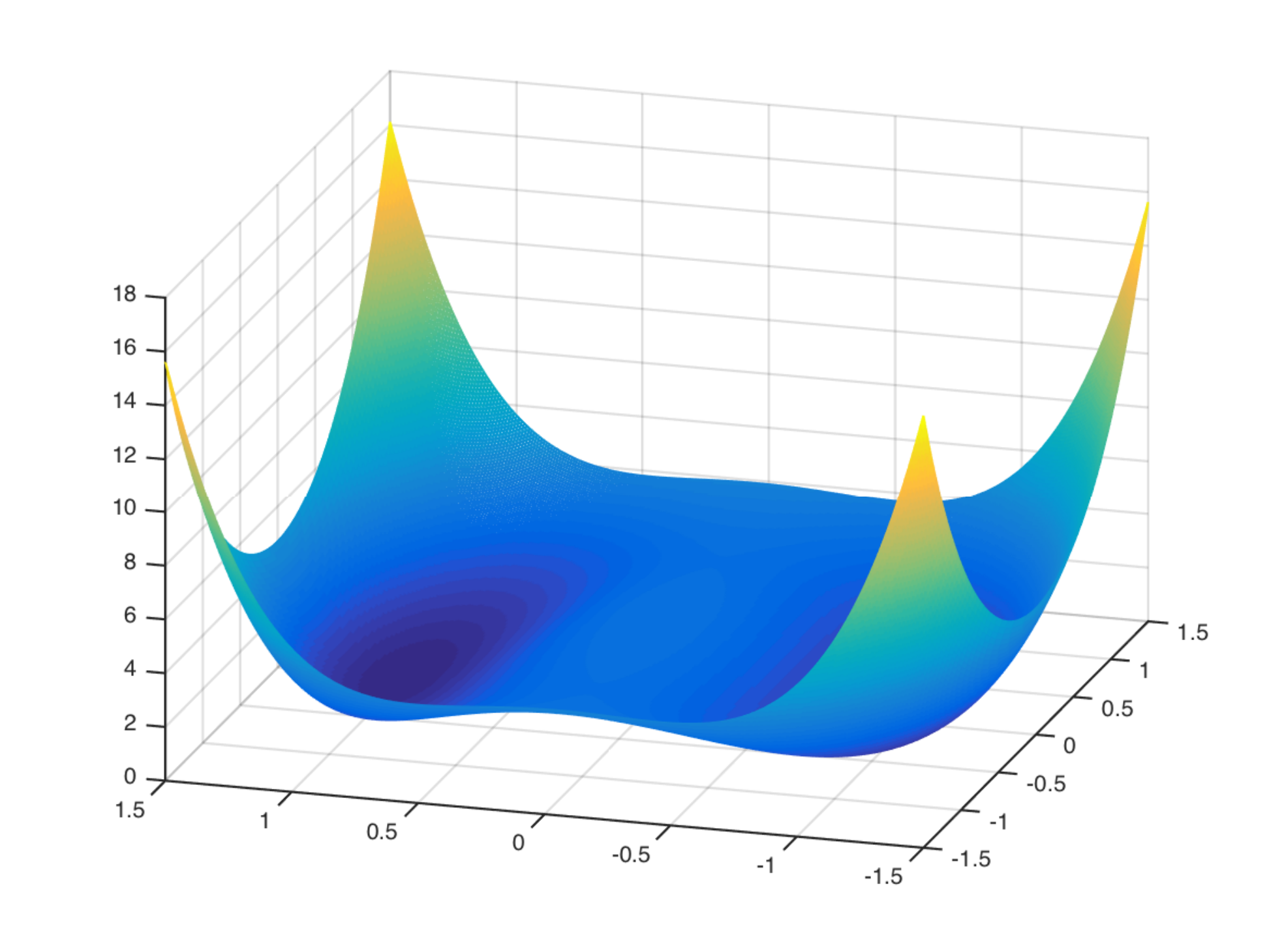}
\caption{(Left) Function landscape of learning sparsifying complete basis via~\eqref{eq:dl_form} in $\R^3$. (Right) Function landscape of generalized phase retrieval via~\eqref{eq:pr_lsq}, assuming the target signal $\mb x$ is real in $\R^2$. In each case, note the equivalent global minimizers and the ridable saddles.  }
\end{SCfigure}

\begin{itemize}
\item \textbf{The Eigenvector Problem.} For a symmetric matrix $\mb A \in \R^{n \times n}$, the classic eigenvector problem is 
\begin{align}
\maxi_{\mb x \in \R^n} \; \mb x^\top \mb A \mb x \quad \st \quad \|\mb x\| = 1. 
\end{align}
Here the manifold is the sphere $\bb S^{n-1}$. It can be easily shown that (see, e.g., Section 4.6 of~\cite{absil2009}) the set of critical points to the problem is exactly the set of eigenvectors to $\mb A$. Moreover, suppose $\lambda_1 > \lambda_2 \ge \dots \lambda_{n-1} > \lambda_n$, with the corresponding eigenvectors $\mb v_1, \dots, \mb v_n$. Then, the only global maximizers are $\pm \mb v_1$, the only global minimizers are $\pm \mb v_n$, and all the intermediate eigenvectors and their negatives are ridable saddle points.\footnote{One can state a more general version of the results, allowing multiplicity of maximum eigenvalues.} Quantitatively, one can show that the function is $(c(\lambda_{n-1} - \lambda_n), c(\lambda_{n-1} - \lambda_n)/\lambda_1, c(\lambda_{n-1} - \lambda_n), 2c(\lambda_{n-1} - \lambda_n)/\lambda_1)$-ridable over $\bb S^{n-1}$ for a certain absolute constant $c > 0$. 

\item \textbf{Complete Dictionary Recovery~\cite{sun2015complete_tr}.} Arising in signal processing and machine learning, dictionary learning tries to approximate a given data matrix $\mb Y \in \R^{n \times p}$ as the product of a dictionary $\mb A$ and a sparse coefficient matrix $\mb X$. In recovery setting, assuming $\mb Y = \mb A_0 \mb X_0$ with $\mb A_0$ square and invertible, $\mb Y$ and $\mb X_0$ have the same row space. Under appropriate model on $\mb X_0$, it makes sense to recover one row of $\mb X_0$ each time by finding the sparsest direction\footnote{The absolute scale is not recoverable. } in $\mathrm{row}(\mb Y)$ by solving the optimization:
\begin{align*}
\mini_{\mb q} \norm{\mb q^\top \mb Y}{0} \quad \st \quad \mb q\neq \mb 0, 
\end{align*}
which can be relaxed as  
\begin{align} \label{eq:dl_form}
\mini f(\mb q) \doteq \frac{1}{p}\sum_{k=1}^p h(\mb q^\top \ol{\mb y}_k) \quad \st \quad \norm{\mb q}{2} = 1 \quad [\text{i.e.,}\; \mb q \in \bb S^{n-1}]. 
\end{align}
Here $h(\cdot)$ is a smooth approximation to the $\abs{\cdot}$ function and $\ol{\mb y}_k$ the $k$-th column of $\ol{\mb Y}$, a proxy of $\mb Y$. The manifold $\mc M$ is $\bb S^{n-1}$ here. \cite{sun2015complete_tr} (Theorem 2.3 and Corollary 2.4) showed that when $h(\cdot) = \mu \log \cosh(\cdot/\mu)$ and $p$ is reasonably large, those $\mb q$'s that help recover rows of $\mb X_0$ are \emph{the only local minimizers} of $f$ over $\bb S^{n-1}$.\footnote{These local minimizers are all global when $p \to \infty$. For finite $p$ that is large enough, these local minimizers assume very close values, and each of them produces a close approximation to a row of $\mb X_0$. } Moreover, these exists a positive constant $c$ such that $f$ is $(c\theta, c\theta, c\theta/\mu, \sqrt{2}\mu/7)$-ridable over $\bb S^{n-1}$, where $\theta$ controls the sparsity level of $\mb X_0$. 
\item \textbf{(Generalized) Phase Retrieval~\cite{sun2015geometric}.} For complex signal $\mb x \in \Cp^n$, generalized phase retrieval (GPR) tries to recover $\mb x$ from the nonlinear measurements of the form $y_k = \abs{\mb a_k^* \mb x}$, for $k = 1, \dots, m$. This task has occupied the central place in imaging systems for scientific discovery~\cite{shechtman2015phase}. Assuming i.i.d. Gaussian measurement noise, a natural formulation for GPR is 
\begin{align} \label{eq:pr_lsq}
\mini_{\mb z \in \Cp^n} f(\mb z) \doteq  \frac{1}{4m} \sum_{k=1}^m (y_k^2 - \abs{\mb a_k^* \mb z}^2)^2. 
\end{align}
The manifold $\mc M$ here is $\Cp^n$. It is obvious that for all $\mb z$, $f(\mb z)$ has the same value as $f(\mb z \e^{\im \theta})$ for any $\theta \in [0, 2\pi)$. \cite{sun2015geometric} showed when $m \ge \Omega(n \log^3 n)$, $\set{\mb x \e^{\im \theta}}$ are the only local minimizers, and also global minimizers (as $f \ge 0$). Moreover, modulo the trivial equivalence discussed above, the function $f$ is $(c, c/(n \log m), c, c/(n \log m))$-ridable for a certain absolute constant $c$, assuming $\norm{\mb x}{} = 1$. 

\item \textbf{Independent Component Analysis (ICA) and Orthogonal Tensor Decomposition~\cite{ge2015escaping}.} Typical setting of ICA asks for a linear transformation $\mb A$ for a given data matrix $\mb Y$, such that rows of $\mb A \mb Y$ achieve maximal statistical independence. Tensor decomposition generalizes spectral decomposition of matrices. Here we focus on \emph{orthogonally decomposable} $d$-th order tensors $\mc T$ which can be represented as 
\begin{align*}
\mc T = \sum_{i=1}^n \mb a_i^{\otimes d}, \quad \mb a_i^\top \mb a_j = \delta_{ij} \; \forall\; i, j,      (\mb a_i \in \R^n\; \forall\; i)
\end{align*}
where $\otimes$ generalizes the usual outer product of vectors. Tensor decomposition refers to finding (up to sign and permutation) the components $\mb a_i$'s given $\mc T$. With appropriate processing and up to small perturbation, ICA is showed to be equivalent to decomposition of a certain form of $4$-th order orthogonally decomposable tensors~\cite{frieze1996learning, arora2012provable}. Specifically, \cite{ge2015escaping} showed (Section C.1.) \footnote{\cite{ge2015escaping} has not used the manifold language as we use here, but resorted to Lagrange multiplier and optimality of the Lagrangian function. For the two decomposition formulations we discussed here, one can verify that the gradient and Hessian they defined are exactly the Riemannian gradient and Hessian of the respective manifolds.  } the minimization problem 
\begin{align*}
\mini f(\mb u) \doteq -\mc T(\mb u, \mb u, \mb u, \mb u) = -\sum_{i=1}^n (\mb a_i^\top \mb u)^4 \quad \st \quad \norm{\mb u}{2} = 1
\end{align*}
has $\pm \mb a_i$'s as its only minimizers and the function $f$ is $(7/n, 1/\mathrm{poly}(n), 3, 1/\mathrm{poly}(n))$-ridable over $\bb S^{n-1}$. Once one of the component is obtained, one can apply deflation to obtain the others. One alternative that tends to make the process more noise-stable is trying to recover all the components in one shot. To this end, \cite{ge2015escaping} proposed to solve 
\begin{align*}
& \mini g(\mb u_1, \dots, \mb u_r) \doteq \sum_{i \ne j} \mc T(\mb u_i, \mb u_i, \mb u_j, \mb u_j) = \sum_{i \ne j} \sum_{k=1}^n (\mb a_k^\top \mb u_i)^2 (\mb a_k^\top \mb u_j)^2, \\
& \st\; \norm{\mb u_i}{} = 1\; \forall i \in [n]. 
\end{align*}
The object $\set{\mb U \in \R^{n \times n}: \norm{\mb u_i}{} = 1\; \forall i}$ is called the \emph{oblique manifold}, which is a product space of multiple spheres. \cite{ge2015escaping} showed all local minimizers of $g$ are equivalent (i.e., signed permuted) copies of $[\mb a_1, \dots, \mb a_n]$. Moreover, $g$ is $(1/\mathrm{poly}(n), 1/\mathrm{poly}(n), 1, 1/\mathrm{poly}(n))$-ridable. 

\item \textbf{Phase Synchronization and Community Detection~\cite{boumal2016nonconvex,bandeira2016low}.} 
Phase synchronization concerns recovery of unit-modulus complex scalars from their relative phases. More precisely, recovering an unknown vector $\mb z \in \Cp_1^n$ with 
\begin{align*}
\Cp_1^n \doteq \set{\mb z \in \Cp^n: \abs{z_1} = \dots = \abs{z_n} = 1},  
\end{align*}
from noisy measurements of the form $C_{ij} = z_i \ol{z_j} + \Delta_{ij}$. The problem is interesting when the noise is nonzero yet controlled, which demands robust solution schemes. Turning to the optimization approach, a natural formulation (if one assumes a Gaussian noise model) is 
\begin{align*}
\mini_{\mb x \in \Cp_1^n} \; \norm{\mb x \mb x^* - \mb C}{F}^2, 
\end{align*}
where we have collected $C_{{ij}}$ into a matrix $\mb C$. Assuming the noise is symmetric (i.e., $\Delta_{ij} = \Delta_{ji}$), the above formulation is equivalent to 
\begin{align} \label{eq:morex_phase_sync}
\mini_{\mb x \in \Cp_1^n} \; -\mb x^* \mb C \mb x. 
\end{align}

Interestingly, for the phase synchronization model, i.e., $\mb C = \mb z \mb z^* + \mb \Delta$ with Hermitian noise matrix $\mb \Delta$, \cite{boumal2016nonconvex} recently showed that (Theorem 4) when the noise $\mb \Delta$ is bounded in mild sense, 
\begin{quote}
\centering
second-order necessary condition for optimality is also sufficient. 
\end{quote}
Particularly, this holds w.h.p. when the noise is i.i.d. complex Gaussians with small variance (Lemma 5). To understand the above statement, recall that second-order necessary condition asks for vanishing gradient and negative semidefinite Hessian at a point. The above statement asserts that such condition is sufficient to guarantee global optimality. In other words, at any critical points other than these verifying the condition have indefinite Hessians. Thus, \cite{boumal2016nonconvex} has effectively showed that when $\mb \Delta$ is appropriately bounded, 
\begin{quote}
\centering 
the function $-\mb x^* \mb C \mb x$ over $\Cp_1^n$ is a ``qualitative'' $\mc X$ function. \footnote{Strictly speaking, our definition of $\mc X$ functions requires the function to be locally \emph{strongly} convex around the local/global minimizers, while the Hessian being positive semidefinite is weaker than that. No matter whether their result can be strengthened in this respect, we note that we impose the strong convexity assumption instead of just convexity is for the sake of deriving concrete convergence rates for optimization algorithms. One can relax the requirement when talking of the qualitative aspect of the structure. Similar comment applies to the ensuring discussion of the real version also. }
\end{quote}

The real counterpart of phase synchronization is called \emph{synchronization over $\Z^2$}, i.e., $\mb z \in \set{1, -1}^n$. In this case, an analogous formulation to~\eqref{eq:morex_phase_sync} appears to be a hard combinatorial problem (think of \texttt{MAX-CUT}) in theory, and also not be friendly for numerical computation (the domain is discrete). Interestingly, \cite{bandeira2016low} showed certain \emph{nonconvex} relaxation has a benign geometric structure. Specifically, applying the usual SDP lifting idea leads to 
\begin{align*}
\mini_{\mb X \in \R^{n \times n}}\; -\innerprod{\mb X}{\mb C} \quad X_{ii} = 1, \forall\; i, \;\; \mb X \succeq \mb 0, \; \; \rank(\mb X) = 1. 
\end{align*}
Dropping the rank constraint results in a convex program (SDP), which is expensive to solve for large $n$. The Burer-Monteiro factorization approach~\cite{burer2003nonlinear,burer2005local} suggests substituting $\mb X = \mb W \mb W^\top$ for $\mb W \in \R^{n \times p}$ for $1 \le p \ll n$ such that the above relaxation is reformulated as 
\begin{align} \label{eq:morex_ncvx_relax_z2}
\mini_{\mb W \in \R^{n \times p}} \; - \trace\paren{\mb W^\top \mb C \mb W} \quad \norm{\mb w_i}{} = 1\; \forall\; i. 
\end{align}
Classic results~\cite{shapiro1982rank,barvinok1995problems,pataki1998rank} on this says \eqref{eq:morex_ncvx_relax_z2} has the same optimal value as the SDP relaxation when $p$ is large enough ($p \sim \Theta(\sqrt{n})$). Moreover, when $p$ is set to be this scale, rank-deficient local optimizers are also global~\cite{burer2005local}. Surprisingly, \cite{bandeira2016low} showed (Theorem 4) that even $p=2$, for the $\Z^2$ synchronization problem with small noise (i.e., small $\mb \Delta$), formulation~\eqref{eq:morex_ncvx_relax_z2} obeys 
\begin{quote}
all points verifying the second-order necessary condition are global optimizes.  
\end{quote}
By analogous argument as for the complex case, this implies: 
\begin{quote}
the function $- \trace\paren{\mb W^\top \mb C \mb W}$ over the oblique manifold $\set{\mb W \in \R^{n \times 2}: \norm{\mb w_i}{}=1 \; \forall i \in [n]}$ is a qualitative $\mc X$ function. 
\end{quote}
A similar result was derived in~\cite{bandeira2016low} for the two-block community detection problem based on the stochastic block model (Theorem 6).\footnote{Both~\cite{bandeira2016low} and~\cite{montanari2016grothendieck} also contain results that characterize local optimizers in terms of their correlation with the optimizer under less stringent/general conditions. }



\end{itemize}

\section{Second-order Trust-region Method and Proof of Convergence}
The intuition that second-order information can help escape ridable saddles from the very start suggests a second-order method. We describe a \emph{second-order trust-region algorithm on manifolds}~\cite{absil2007trust, absil2009} for this purpose. 

For the generic problem~\eqref{eq:manifold_min}, we start from any feasible $\mb x^{(0)} \in \mc M$, and form a sequence of iterates $\mb x^{(1)}, \mb x^{(2)}, \dots \in \mc M$ as follows. For the current iterate $\mb x^{(k)}$, we consider the quadratic approximation 
\begin{align} \label{eq:tr_approx}
\wh{f}(\mb \delta; \mb x^{(k)}) \doteq f(\mb x^{(k)}) + \innerprod{\mb \delta}{\grad f(\mb x^{(k)})} + \frac{1}{2} \innerprod{\Hess f(\mb x^{(k)})[\mb \delta]}{\mb \delta}
\end{align}
which is defined for all $\mb \delta \in T_{\mb x^{(k)}} \mc M$. The next iterate is determined by minimizing the quadratic approximation within a small radius $\Delta$ (i.e., the trust region) of $\mb x^{(k)}$, i.e., 
\begin{align} \label{eq:rie_tr_subprob}
\mb \delta^{(k+1)} \doteq \mathop{\arg\min}_{\mb \delta \in T_{\mb x^{(k)}} \mc M, \norm{\mb \delta}{2} \le \Delta} \widehat{f}\paren{\mb \delta; \mb x^{(k)}}, 
\end{align}
which is called the Riemannian trust-region subproblem. The vector $\mb x^{(k)} + \mb \delta^{(k+1)} $ is generally not a point on $\mc M$. One then performs a \emph{retraction} step $R_{\mb x^{(k)}}$ that pulls the vector back to the manifold, resulting in the update formula
\begin{align*}
\mb x^{(k+1)} = R_{\mb x^{(k)}} (\mb x^{(k)} + \mb \delta^{(k+1)}). 
\end{align*}
Most manifolds of practical interest are embedded submanifolds of $\R^{m \times n}$ and the tangent space is a subspace of $\R^{m \times n}$. For an $\mb x^{(k)} \in \mc M$ and an orthonormal basis $\mb U$ for $T_{\mb x^{(k)}} \mc M$, one can solve~\eqref{eq:rie_tr_subprob} by solving the recast Euclidean trust-region subproblem 
\begin{align}
\mb \xi^{(k+1)} \doteq \mathop{\arg\min}_{\norm{\mb \xi}{} \le \Delta} \wh{f}(\mb U \mb \xi; \mb x^{(k)}), 
\end{align}
for which efficient numerical algorithms exist~\cite{more1983computing, conn2000trust, fortin2004trust, hazan2014linear}. 
\begin{wrapfigure}{R}{0.25\textwidth}
\centering 
\includegraphics[width=0.9\linewidth]{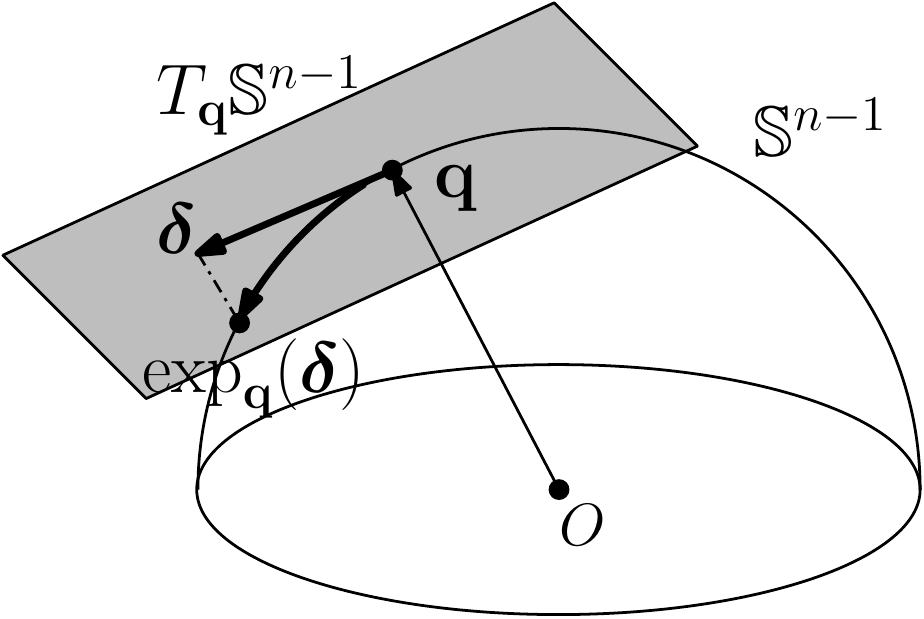}
\caption{Illustrations of the tangent space $T_{\mb q}\bb S^{n-1}$ and exponential map $\exp_{\mb q}\paren{\mb \delta}$ defined on the sphere $\bb S^{n-1}$.}
 \label{fig:exp-map}
\end{wrapfigure}
Design choice of the retraction is often problem-specific, ranging from the classical exponential map to the Euclidean projection that works for many matrix manifolds~\cite{absil2012projection}. 

To show the trust-region algorithm converges to a global minimizer, we assume $\Delta$ is small enough such that approximation error of~\eqref{eq:tr_approx} to $f$ is ``negligible'' locally. Each step around a negative-curvature or strong-gradient point decreases the objective by a certain amount. Indeed, it is clear there is always one descent direction in such cases. Thus, the trust-region step will approximately follow one descent direction and decrease the function value. When the iterate sequence moves into a strongly convex region around a global minimizer, a step is either constrained such that it also deceases the objective by an amount, or unconstrained, which is a good indicator that the target minimizer is within a radius $\Delta$. In the latter case, the algorithm behaves like the classical Newton method and quadratic sequence convergence can be shown. 

Quantitative convergence proof demands knowledge of the ridability parameters, smoothness parameters of the objective, and elements of Riemannian geometry. We refer the reader to~\cite{sun2015complete_tr,sun2015geometric} for practical examples of convergence analyses. 

\section{Discussion}
Recently, there is a surge of interest in understanding nonconvex heuristics for practical problems~\cite{keshavan2010matrix, jain2013low, hardt2014understanding, hardt2014fast, netrapalli2014non, jain2014fast, sun2014guaranteed, zheng2015convergent, tu2015low, chen2015fast, netrapalli2013phase, candes2015phase, chen2015solving, white2015local, jain2014provable, anandkumar2014guaranteed, anandkumar2014analyzing, anandkumar2015tensor, yi2013alternating, sedghi2014provable, lee2013near, qu2014finding, lee2013near, agarwal2013learning, agarwal2013exact, arora2013new, arora2015simple, arora2014more, jain2015computing}. Majority of the work start from clever initializations, and then proceed with analysis of local convergence. In comparison, it is clear that for $\mc X$ functions, second-order trust-region algorithms with any initialization guarantee to retrieve one target minimizer. Identifying $\mc X$ functions has involved intensive technical work~\cite{sun2015complete_tr,sun2015geometric,ge2015escaping}. It is interesting to see if streamlined toolkits can be developed, say via operational rules or unified potential functions. This would facilitate study of other practical problems, such as the deep networks of which saddle points are believed to be prevalent and constitute significant computational bottleneck~\cite{pascanu2014saddle, dauphin2014identifying, choromanska2014loss}. To match heuristics computationally, more practical algorithms other than the second-order trust-region methods are needed. Practical trust-region solvers with saddle-escaping capability may be possible for structured problems~\cite{boumal2014manopt,sun2015geometric}. Moreover, simulations with several practical problems suggest gradient-style algorithms with random initializations succeed. \cite{ge2015escaping,lee2016gradient} are recent endeavors towards this direction.

{\small
\bibliographystyle{amsalpha}
\bibliography{ncvx,dl_tit,pr,thesis,talk}

\newcommand{\etalchar}[1]{$^{#1}$}
\providecommand{\bysame}{\leavevmode\hbox to3em{\hrulefill}\thinspace}
\providecommand{\MR}{\relax\ifhmode\unskip\space\fi MR }
\providecommand{\MRhref}[2]{%
  \href{http://www.ams.org/mathscinet-getitem?mr=#1}{#2}
}
\providecommand{\href}[2]{#2}
\begin{thebibliography}{AGMM15}

\bibitem[AAJ{\etalchar{+}}13]{agarwal2013learning}
Alekh Agarwal, Animashree Anandkumar, Prateek Jain, Praneeth Netrapalli, and
  Rashish Tandon, \emph{Learning sparsely used overcomplete dictionaries via
  alternating minimization}, arXiv preprint arXiv:1310.7991 (2013).

\bibitem[AAN13]{agarwal2013exact}
Alekh Agarwal, Animashree Anandkumar, and Praneeth Netrapalli, \emph{Exact
  recovery of sparsely used overcomplete dictionaries}, arXiv preprint
  arXiv:1309.1952 (2013).

\bibitem[ABG07]{absil2007trust}
Pierre-Antoine. Absil, Christopher~G. Baker, and Kyle~A. Gallivan,
  \emph{Trust-region methods on {R}iemannian manifolds}, Foundations of
  Computational Mathematics \textbf{7} (2007), no.~3, 303--330.

\bibitem[ABGM14]{arora2014more}
Sanjeev Arora, Aditya Bhaskara, Rong Ge, and Tengyu Ma, \emph{More algorithms
  for provable dictionary learning}, arXiv preprint arXiv:1401.0579 (2014).

\bibitem[AG16]{anandkumar2016efficient}
Anima Anandkumar and Rong Ge, \emph{Efficient approaches for escaping higher
  order saddle points in non-convex optimization}, arXiv preprint
  arXiv:1602.05908 (2016).

\bibitem[AGJ14a]{anandkumar2014analyzing}
Animashree Anandkumar, Rong Ge, and Majid Janzamin, \emph{Analyzing tensor
  power method dynamics: Applications to learning overcomplete latent variable
  models}, arXiv preprint arXiv:1411.1488 (2014).

\bibitem[AGJ14b]{anandkumar2014guaranteed}
\bysame, \emph{Guaranteed non-orthogonal tensor decomposition via alternating
  rank-1 updates}, arXiv preprint arXiv:1402.5180 (2014).

\bibitem[AGM13]{arora2013new}
Sanjeev Arora, Rong Ge, and Ankur Moitra, \emph{New algorithms for learning
  incoherent and overcomplete dictionaries}, arXiv preprint arXiv:1308.6273
  (2013).

\bibitem[AGMM15]{arora2015simple}
Sanjeev Arora, Rong Ge, Tengyu Ma, and Ankur Moitra, \emph{Simple, efficient,
  and neural algorithms for sparse coding}, arXiv preprint arXiv:1503.00778
  (2015).

\bibitem[AGMS12]{arora2012provable}
Sanjeev Arora, Rong Ge, Ankur Moitra, and Sushant Sachdeva, \emph{Provable
  {ICA} with unknown gaussian noise, with implications for gaussian mixtures
  and autoencoders}, Advances in Neural Information Processing Systems, 2012,
  pp.~2375--2383.

\bibitem[AJSN15]{anandkumar2015tensor}
Animashree Anandkumar, Prateek Jain, Yang Shi, and Uma~Naresh Niranjan,
  \emph{Tensor vs matrix methods: Robust tensor decomposition under block
  sparse perturbations}, arXiv preprint arXiv:1510.04747 (2015).

\bibitem[AM12]{absil2012projection}
P.-A. Absil and J{\'e}r{\^o}me Malick, \emph{Projection-like retractions on
  matrix manifolds}, SIAM Journal on Optimization \textbf{22} (2012), no.~1,
  135--158.

\bibitem[AMS09]{absil2009}
Pierre-Antoine. Absil, Robert Mahoney, and Rodolphe Sepulchre,
  \emph{Optimization algorithms on matrix manifolds}, Princeton University
  Press, 2009.

\bibitem[Bar95]{barvinok1995problems}
Alexander~I. Barvinok, \emph{Problems of distance geometry and convex
  properties of quadratic maps}, Discrete \& Computational Geometry \textbf{13}
  (1995), no.~2, 189--202.

\bibitem[BBV16]{bandeira2016low}
Afonso~S Bandeira, Nicolas Boumal, and Vladislav Voroninski, \emph{On the
  low-rank approach for semidefinite programs arising in synchronization and
  community detection}, arXiv preprint arXiv:1602.04426 (2016).

\bibitem[Ber99]{bertsekas1999nonlinear}
Dimitri~P. Bertsekas, \emph{Nonlinear programming}, Athena scientific, 1999.

\bibitem[BM03]{burer2003nonlinear}
Samuel Burer and Renato~D.C. Monteiro, \emph{A nonlinear programming algorithm
  for solving semidefinite programs via low-rank factorization}, Mathematical
  Programming \textbf{95} (2003), no.~2, 329--357.

\bibitem[BM05]{burer2005local}
Samuel Burer and Renato D.~C. Monteiro, \emph{Local minima and convergence in
  low-rank semidefinite programming}, Mathematical Programming \textbf{103}
  (2005), no.~3, 427--444.

\bibitem[BMAS14]{boumal2014manopt}
Nicolas Boumal, Bamdev Mishra, P.-A. Absil, and Rodolphe Sepulchre,
  \emph{{M}anopt, a {M}atlab toolbox for optimization on manifolds}, Journal of
  Machine Learning Research \textbf{15} (2014), 1455--1459.

\bibitem[Bou16]{boumal2016nonconvex}
Nicolas Boumal, \emph{Nonconvex phase synchronization}, arXiv preprint
  arXiv:1601.06114 (2016).

\bibitem[CC15]{chen2015solving}
Yuxin Chen and Emmanuel~J. Cand{\`e}s, \emph{Solving random quadratic systems
  of equations is nearly as easy as solving linear systems}, arXiv preprint
  arXiv:1505.05114 (2015).

\bibitem[CGT00]{conn2000trust}
Andrew~R. Conn, Nicholas I.~M. Gould, and Philippe~L. Toint, \emph{Trust-region
  methods}, Society for Industrial and Applied Mathematics, Philadelphia, PA,
  USA, 2000.

\bibitem[CHM{\etalchar{+}}14]{choromanska2014loss}
Anna Choromanska, Mikael Henaff, Michael Mathieu, G{\'e}rard~Ben Arous, and
  Yann LeCun, \emph{The loss surface of multilayer networks}, arXiv preprint
  arXiv:1412.0233 (2014).

\bibitem[CLS15]{candes2015phase}
Emmanuel~J. Cand{\`e}s, Xiaodong Li, and Mahdi Soltanolkotabi, \emph{Phase
  retrieval via wirtinger flow: Theory and algorithms}, Information Theory,
  IEEE Transactions on \textbf{61} (2015), no.~4, 1985--2007.

\bibitem[CW15]{chen2015fast}
Yudong Chen and Martin~J. Wainwright, \emph{Fast low-rank estimation by
  projected gradient descent: General statistical and algorithmic guarantees},
  arXiv preprint arXiv:1509.03025 (2015).

\bibitem[DPG{\etalchar{+}}14]{dauphin2014identifying}
Yann~N. Dauphin, Razvan Pascanu, Caglar Gulcehre, Kyunghyun Cho, Surya Ganguli,
  and Yoshua Bengio, \emph{Identifying and attacking the saddle point problem
  in high-dimensional non-convex optimization}, Advances in Neural Information
  Processing Systems, 2014, pp.~2933--2941.

\bibitem[FJK96]{frieze1996learning}
Alan Frieze, Mark Jerrum, and Ravi Kannan, \emph{Learning linear
  transformations}, focs, IEEE, 1996, p.~359.

\bibitem[FW04]{fortin2004trust}
Charles Fortin and Henry Wolkowicz, \emph{The trust region subproblem and
  semidefinite programming}, Optimization methods and software \textbf{19}
  (2004), no.~1, 41--67.

\bibitem[GHJY15]{ge2015escaping}
Rong Ge, Furong Huang, Chi Jin, and Yang Yuan, \emph{Escaping from saddle
  points---online stochastic gradient for tensor decomposition}, Proceedings of
  The 28th Conference on Learning Theory, 2015, pp.~797--842.

\bibitem[Har14]{hardt2014understanding}
Moritz Hardt, \emph{Understanding alternating minimization for matrix
  completion}, Foundations of Computer Science (FOCS), 2014 IEEE 55th Annual
  Symposium on, IEEE, 2014, pp.~651--660.

\bibitem[HK14]{hazan2014linear}
Elad Hazan and Tomer Koren, \emph{A linear-time algorithm for trust region
  problems}, arXiv preprint arXiv:1401.6757 (2014).

\bibitem[HL13]{hillar2013most}
Christopher~J. Hillar and Lek-Heng Lim, \emph{Most tensor problems are
  {NP}-hard}, Journal of the ACM (JACM) \textbf{60} (2013), no.~6, 45.

\bibitem[HW14]{hardt2014fast}
Moritz Hardt and Mary Wootters, \emph{Fast matrix completion without the
  condition number}, Proceedings of The 27th Conference on Learning Theory,
  2014, pp.~638--678.

\bibitem[JJKN15]{jain2015computing}
Prateek Jain, Chi Jin, Sham~M. Kakade, and Praneeth Netrapalli, \emph{Computing
  matrix squareroot via non convex local search}, arXiv preprint
  arXiv:1507.05854 (2015).

\bibitem[JN14]{jain2014fast}
Prateek Jain and Praneeth Netrapalli, \emph{Fast exact matrix completion with
  finite samples}, arXiv preprint arXiv:1411.1087 (2014).

\bibitem[JNS13]{jain2013low}
Prateek Jain, Praneeth Netrapalli, and Sujay Sanghavi, \emph{Low-rank matrix
  completion using alternating minimization}, Proceedings of the forty-fifth
  annual ACM symposium on Theory of Computing, ACM, 2013, pp.~665--674.

\bibitem[JO14]{jain2014provable}
Prateek Jain and Sewoong Oh, \emph{Provable tensor factorization with missing
  data}, Advances in Neural Information Processing Systems, 2014,
  pp.~1431--1439.

\bibitem[KMO10]{keshavan2010matrix}
Raghunandan~H. Keshavan, Andrea Montanari, and Sewoong Oh, \emph{Matrix
  completion from a few entries}, Information Theory, IEEE Transactions on
  \textbf{56} (2010), no.~6, 2980--2998.

\bibitem[LSJR16]{lee2016gradient}
Jason~D Lee, Max Simchowitz, Michael~I Jordan, and Benjamin Recht,
  \emph{Gradient descent converges to minimizers}, arXiv preprint
  arXiv:1602.04915 (2016).

\bibitem[LWB13]{lee2013near}
Kiryung Lee, Yihong Wu, and Yoram Bresler, \emph{Near optimal compressed
  sensing of sparse rank-one matrices via sparse power factorization}, arXiv
  preprint arXiv:1312.0525 (2013).

\bibitem[MK87]{murty1987some}
Katta~G. Murty and Santosh~N. Kabadi, \emph{Some {NP}-complete problems in
  quadratic and nonlinear programming}, Mathematical programming \textbf{39}
  (1987), no.~2, 117--129.

\bibitem[Mon16]{montanari2016grothendieck}
Andrea Montanari, \emph{A grothendieck-type inequality for local maxima}, arXiv
  preprint arXiv:1603.04064 (2016).

\bibitem[MS83]{more1983computing}
Jorge~J. Mor{\'e} and Danny~C. Sorensen, \emph{Computing a trust region step},
  SIAM Journal on Scientific and Statistical Computing \textbf{4} (1983),
  no.~3, 553--572.

\bibitem[NJS13]{netrapalli2013phase}
Praneeth Netrapalli, Prateek Jain, and Sujay Sanghavi, \emph{Phase retrieval
  using alternating minimization}, Advances in Neural Information Processing
  Systems, 2013, pp.~2796--2804.

\bibitem[NNS{\etalchar{+}}14]{netrapalli2014non}
Praneeth Netrapalli, Uma~Naresh. Niranjan, Sujay Sanghavi, Animashree
  Anandkumar, and Prateek Jain, \emph{Non-convex robust {PCA}}, Advances in
  Neural Information Processing Systems, 2014, pp.~1107--1115.

\bibitem[Pat98]{pataki1998rank}
G{\'a}bor Pataki, \emph{On the rank of extreme matrices in semidefinite
  programs and the multiplicity of optimal eigenvalues}, Mathematics of
  operations research \textbf{23} (1998), no.~2, 339--358.

\bibitem[PDGB14]{pascanu2014saddle}
Razvan Pascanu, Yann~N Dauphin, Surya Ganguli, and Yoshua Bengio, \emph{On the
  saddle point problem for non-convex optimization}, arXiv preprint
  arXiv:1405.4604 (2014).

\bibitem[QSW14]{qu2014finding}
Qing Qu, Ju~Sun, and John Wright, \emph{Finding a sparse vector in a subspace:
  Linear sparsity using alternating directions}, Advances in Neural Information
  Processing Systems, 2014, pp.~3401--3409.

\bibitem[RI10]{reich2010nonlinear}
Simeon Reich and Aleksandr~Davidovich Ioffe, \emph{Nonlinear analysis and
  optimization: Optimization}, vol.~2, American Mathematical Soc., 2010.

\bibitem[SA14]{sedghi2014provable}
Hanie Sedghi and Animashree Anandkumar, \emph{Provable tensor methods for
  learning mixtures of classifiers}, arXiv preprint arXiv:1412.3046 (2014).

\bibitem[SEC{\etalchar{+}}15]{shechtman2015phase}
Yoav Shechtman, Yonina~C. Eldar, Oren Cohen, Henry~N. Chapman, Jianwei Miao,
  and Mordechai Segev, \emph{Phase retrieval with application to optical
  imaging: A contemporary overview}, Signal Processing Magazine, IEEE
  \textbf{32} (2015), no.~3, 87--109.

\bibitem[Sha82]{shapiro1982rank}
Alexander Shapiro, \emph{Rank-reducibility of a symmetric matrix and sampling
  theory of minimum trace factor analysis}, Psychometrika \textbf{47} (1982),
  no.~2, 187--199.

\bibitem[SL14]{sun2014guaranteed}
Ruoyu Sun and Zhi-Quan Luo, \emph{Guaranteed matrix completion via non-convex
  factorization}, arXiv preprint arXiv:1411.8003 (2014).

\bibitem[SQW15]{sun2015complete_tr}
Ju~Sun, Qing Qu, and John Wright, \emph{Complete dictionary recovery over the
  sphere}, arXiv preprint arXiv:1504.06785 (2015).

\bibitem[SQW16]{sun2015geometric}
\bysame, \emph{A geometric analysis of phase retreival}, arXiv preprint
  arXiv:1602.06664 (2016).

\bibitem[TBSR15]{tu2015low}
Stephen Tu, Ross Boczar, Mahdi Soltanolkotabi, and Benjamin Recht,
  \emph{Low-rank solutions of linear matrix equations via procrustes flow},
  arXiv preprint arXiv:1507.03566 (2015).

\bibitem[WWS15]{white2015local}
Chris~D. White, Rachel Ward, and Sujay Sanghavi, \emph{The local convexity of
  solving quadratic equations}, arXiv preprint arXiv:1506.07868 (2015).

\bibitem[YCS13]{yi2013alternating}
Xinyang Yi, Constantine Caramanis, and Sujay Sanghavi, \emph{Alternating
  minimization for mixed linear regression}, arXiv preprint arXiv:1310.3745
  (2013).

\bibitem[ZL15]{zheng2015convergent}
Qinqing Zheng and John Lafferty, \emph{A convergent gradient descent algorithm
  for rank minimization and semidefinite programming from random linear
  measurements}, arXiv preprint arXiv:1506.06081 (2015).

\end{thebibliography}
}

\end{document}